%AMS TeX file for the paper
%
%      On Jacobi and continuous Hahn polynomials
%
%by H.T. Koelink
%%%%%%%%%%%%%%%%%%%%%%%%%%%%AMS TeX file%%%%%%%%%%%%%%%%%%%%%%%%%%%%%%%%%
\input amstex.tex
\documentstyle{amsppt}
\magnification=1200
\baselineskip=13pt
\hsize=6.5truein
\vsize=8.9truein
%%%%%%%%%%%%%%%%%%%%%%%%%%Numbering macros%%%%%%%%%%%%%%%%%%%%%%%%%%%%%%%%%
\countdef\sectionno=1
\sectionno=0

\countdef\eqnumber=10
\countdef\theoremno=11
\countdef\countrefno=12

\def\newsection{\global\advance\sectionno by 1
                \global\eqnumber=1
                \global\theoremno=1
                \the\sectionno}

\def\theoremname#1{\the\sectionno.\the\theoremno
                   \xdef#1{{\the\sectionno.\the\theoremno}}
                   \global\advance\theoremno by 1}

\def\eqname#1{\the\sectionno.\the\eqnumber
              \xdef#1{{\the\sectionno.\the\eqnumber}}
              \global\advance\eqnumber by 1}

\global\countrefno=1

\def\refno#1{\xdef#1{{\the\countrefno}}\global\advance\countrefno by 1}

%%%%%%%%%%%%%%%%%%%%%%%%%%%%Abbreviations%%%%%%%%%%%%%%%%%%%%%%%%%%%%%%%%%%%
\def\R{{\Bbb R}}
\def\N{{\Bbb N}}
\def\C{{\Bbb C}}

\def\Zp{{\Bbb Z}_+}
\def\F{{\Cal F}}
\def\a{\alpha}
\def\b{\beta}
\def\g{\gamma}
\def\d{\delta}
\def\G{\Gamma}
\def\th{\tanh}
\def\ch{\cosh}

\def\hf{{1\over 2}}
\def\dx{{d\over{dx}}}

%%%%%%%%%%%%%%%%%%%%%%%%%%%Reference numbering%%%%%%%%%%%%%%%%%%%%%%%%%%%%%%%
\refno{\AppeKF}
\refno{\Aske}
\refno{\Asketwee}
\refno{\AskeWeen}
\refno{\AskeWtwee}
\refno{\AtakS}
\refno{\BadeK}
\refno{\Bail}
\refno{\Bateeen}
\refno{\Batetwee}
\refno{\Batedrie}
\refno{\Batevier}
\refno{\Braf}
\refno{\Carleen}
\refno{\Carltwee}
\refno{\Chih}
\refno{\ErdeTIT}
\refno{\Hard}
\refno{\KalnM}
\refno{\KoekS}
\refno{\Koel}
\refno{\Kooreen}
\refno{\Koortwee}
\refno{\Koordrie}
\refno{\NikiSU}
\refno{\Past}
\refno{\RahmS}
\refno{\Rice}
\refno{\Stieeen}
\refno{\Stietwee}
\refno{\Szeg}
\refno{\Titc}
\refno{\Touc}
\refno{\WhitW}
\refno{\Wils}
\refno{\WymaM}

%%%%%%%%%%%%%%%%%%%%%%%%%%TOP MATTER%%%%%%%%%%%%%%%%%%%%%%%%%%%%%%%%%
\topmatter
\title On Jacobi and continuous Hahn polynomials
\endtitle
\author H.T. Koelink\endauthor
%\affil  Katholieke Universiteit Leuven\endaffil
\address Department of Mathematics, Katholieke Universiteit Leuven,
Celestijnenlaan 200 B, B-3001 Leuven (Heverlee), Belgium\endaddress
\email erik\%twi\%wis\@cc3.KULeuven.ac.be\endemail
\date September 20, 1994\enddate
\thanks Supported by a Fellowship of the Research Council
of the Katholieke Universiteit Leuven. \endthanks
\keywords Jacobi polynomials, continuous Hahn polynomials, Fourier transform
\endkeywords
\subjclass 33C45, 42A38
\endsubjclass
\abstract Jacobi polynomials are mapped onto
the continuous Hahn polynomials by the Fourier transform and the orthogonality
relations for the continuous Hahn polynomials then follow from the
orthogonality relations for the Jacobi polynomials and the Parseval formula.
In a special case this relation dates back to work by Bateman in 1933 and
we follow a part of the historical development for these polynomials.
Some applications of this relation are given.
\endabstract
\endtopmatter
\document

%%%%%%%%%%%%%%%%%%%%%%%%%%%%Beginning of text%%%%%%%%%%%%%%%%%%%%%%%%%%%
\head\newsection. Introduction and history\endhead

In Askey's scheme of hypergeometric orthogonal polynomials we find the Jacobi
polynomials and the continuous Hahn polynomials; see Askey and Wilson
\cite{\AskeWtwee, Appendix} with the correction in \cite{\Aske},
Koekoek and Swarttouw \cite{\KoekS} or Koornwinder \cite{\Koortwee, \S 5}
for information on Askey's scheme. In the hierarchy of
Askey's scheme of hypergeometric orthogonal polynomials the continuous Hahn
polynomials are above the Jacobi polynomials since they have one extra degree
of freedom. In this paper we consider a way of going up in the Askey scheme
from the Jacobi polynomials to the continuous Hahn polynomials by use
of the Fourier transform. This method is a simple extension of some special
cases introduced by Bateman in the 1930's.

In his 1933 paper \cite{\Bateeen} Bateman introduces the polynomial $F_n$
satisfying
$$
F_n\Bigl( \dx\Bigr) \ch^{-1}x = \ch^{-1}x P_n(\th x),
\tag\eqname{\vglBatemanpol}
$$
where $P_n(x)=P_n^{(0,0)}(x)$ is the Legendre polynomial,
cf. \thetag{2.3} for its
definition. This is possible since $\dx \ch^{-1}=-\ch^{-1}\th x$ and
$\dx \th x = 1-\th^2x$.
Bateman derives in \cite{\Bateeen} various properties for the
polynomials $F_n$, such as generating functions, explicit expressions as
hypergeometric series, three-term recurrence relation, difference equations,
integral representations and more, see also \cite{\Batetwee}. One year later
Bateman proves the orthogonality
relations, cf. \cite{\Batedrie}, \cite{\Batetwee},
$$
\int_{-\infty}^\infty{{F_n(ix)F_m(ix)}\over{\ch^2 (\pi x/2)}}\, dx =
\d_{n,m} {4(-1)^n\over{\pi (2n+1)}}
\tag\eqname{\vglorthoBatemanpols}
$$
by a method using the Fourier transform, which
we reproduce for a more general class of polynomials.
The factor $(-1)^n$ on the right hand side of
\thetag{\vglorthoBatemanpols} does
not matter, since it can be removed by rescaling $F_n$ by a factor $i^n$.
This factor is necessary in order to make the polynomial $F_n$ real-valued
for imaginary argument.
The orthogonality relation \thetag{\vglorthoBatemanpols} is also derived by
Hardy \cite{\Hard, \S 8} as an example of a general approach to some
orthogonal polynomials using the Mellin transform, which is equivalent to
Bateman's proof of \thetag{\vglorthoBatemanpols}; see also remark 3.2(ii).

The Bateman polynomial $F_n$ has been generalised by Pasternack \cite{\Past}
in 1939. He defines the polynomial $F_n^m$ by
$$
F_n^m\Bigl( \dx\Bigr) \ch^{-m-1}x = \ch^{-m-1}x P_n(\th x),
\tag\eqname{\vglPasternackpol}
$$
for $m\in\C\backslash \{ -1\}$, which reduces to Bateman's polynomial $F_n$
in case $m=0$. Here we use $\dx \ch^{-1-m}x = -(m+1)\ch^{-1-m}x\th x$.
The case $m=1$ already occurs in Bateman's paper
\cite{\Batetwee, \S 4}, see also \cite{\Batevier}.

For the polynomials $F_n^m$ Pasternack derives explicit expressions in terms
of hypergeometric series, generating functions, three-term recurrence
relation, difference equations and integral representations much along
the same lines as in \cite{\Bateeen}, but he does not prove orthogonality
relations for $F_n^m$. In particular, Pasternack proves,
cf. \cite{\Past, (10.2), (10.5)},
$$
\int_{-\infty}^\infty \Bigl\vert \G\bigl(\hf(m+1+z)\bigr)\Bigr\vert^4
F_n^m(iz) F_p^m(-iz)\, dz
= C \int_{-1}^1 (1-x^2)^m P_n(x)P_p(x)\, dx
$$
for some explicit constant $C$. So if Pasternack would have replaced the
Legendre polynomial $P_n$ in \thetag{\vglPasternackpol} by the Gegenbauer
polynomial $P_n^{(m,m)}$, i.e. the polynomials orthogonal on $[-1,1]$ with
respect to the weight function $(1-x^2)^m$, cf. \thetag{2.2}, he would have
obtained a one-parameter subclass of the two-parameter
continuous symmetric Hahn polynomials introduced by Askey and Wilson
\cite{\AskeWeen} in 1982.

Bateman obtains in \cite{\Batevier, (3.3)} the orthogonality relations
$$
\int_{-\infty}^\infty {{F^m_n(ix)F_p^{-m}(ix)}\over
{\ch (\pi x) + \cos (m\pi)}}\, dx =
\d_{n,p} {{2(-1)^n}\over{\pi (2n+1)}}
{{m\pi}\over{\sin\pi m}}, \qquad -1 < m <1.
\tag\eqname{\vglbiorthoPasternackpols}
$$
for Pasternack's polynomials, which reduce to
\thetag{\vglorthoBatemanpols} for $m=0$. The right hand side of
\cite{\Batevier, (3.3)} is not correct.
Although explicit expressions for $F_n^m$, and hence for the leading
coefficients of these polynomials, were known to Bateman at that time he
does not rewrite the orthogonality relation
\thetag{\vglbiorthoPasternackpols} as the orthogonality relation
\thetag{1.5} for the Pasternack polynomials.

In 1956 Touchard \cite{\Touc, \S\S 13, 14} considers orthogonal polynomials
associated with the Bernoulli numbers. He derives a three-term
recurrence relation, orthogonality relations and determines the value at
$-{1\over 2}$. In the following paper in the Canadian Journal of
Mathematics Wyman and Moser \cite{\WymaM} give an explicit expression
for these polynomials in terms of a hypergeometric ${}_4F_3$-series. A
year later Brafman \cite{\Braf} gives another expression for these
polynomials and derives generating functions for these polynomials.
But it is Carlitz
\cite{\Carleen} who notes that the polynomials introduced by Touchard
are the same as Bateman's polynomials $F_n$ defined by
\thetag{\vglBatemanpol}. In 1959 Carlitz \cite{\Carltwee} replaces
the Bernoulli numbers by certain numbers involving the
Bernoulli polynomials at a point $\lambda$, where the case $\lambda=0$
corresponds to the Bernoulli numbers,
see also Chihara \cite{\Chih, Ch.~VI, \S 8}.
Carlitz shows that the corresponding orthogonal
polynomials are Pasternack's polynomials \thetag{\vglPasternackpol}, for
which he gives the orthogonality relations, see also \cite{\AskeWeen},
$$
\int_{-\infty}^\infty {{F_n^m (ix) F_p^m(ix)}\over
{\cos(\pi m) + \ch (\pi x)}}\, dx = \d_{n,p}
{{(-1)^n}\over{2n+1}} {2\over{\pi}} {{(1-m)_n}\over{(1+m)_n}}
{{m\pi}\over{\sin\pi m}}, \qquad -1 < m <1,
\tag\eqname{\vglorthoPasternackpols}
$$
where we use the notation of \thetag{2.4}.
The same remark as for \thetag{\vglorthoBatemanpols} applies here.
The case $m=0$ of \thetag{\vglorthoPasternackpols} is
\thetag{\vglorthoBatemanpols}. From proposition 3.1 we see that
\thetag{\vglorthoPasternackpols} also yields orthogonal polynomials, after
a suitable renormalisation, for $m\in i\R$. From either
\thetag{\vglbiorthoPasternackpols} or from the fact that the weight
function in \thetag{\vglorthoPasternackpols} is even in $m$ we see that
$F_n^m$ is a multiple of $F_n^{-m}$.
By comparing leading coefficients, cf. \S 2, we obtain
$(1+m)_nF_n^m(x)=(1-m)_nF_n^{-m}(x)$, see \S 2 for the notation.
The case $m=\hf$ of \thetag{\vglorthoPasternackpols} has been given by
Hardy \cite{\Hard, \S\S 4-7} in 1940.

Carlitz explicitly calculates the moments corresponding to the orthogonality
measure \thetag{\vglorthoPasternackpols} in terms of the Bernoulli
polynomials, cf. \cite{\Carltwee, \S 6}.
This result has already been obtained by Stieltjes \cite{\Stieeen, \S 5}
in 1890 by developing $\psi(x+b)-\psi(x+1-b)$, $\psi(x)$ denoting the
logarithmic derivative of the gamma function $\Gamma(x)$, in powers of
$x^{-1}$ and in terms of a continued fraction. This means that Stieltjes
gives the three-term recurrence relation for the orthogonal polynomials
for which the moments of the orthogonality measure are given in terms of
the Bernoulli polynomials. Stieltjes includes this example in his
famous memoir ``Recherches sur les fraction continues'', in which he also
gives the corresponding integral representation, cf. \cite{\Stietwee, \S 86}.

In 1982 Askey and Wilson \cite{\AskeWeen}
introduce orthogonal
polynomials, which generalise the orthogonal polynomials introduced by
Bateman, Pasternack, Touchard, Hardy and Carlitz. These polynomials are
orthogonal on $\R$ with respect to the measure
$|\G(\a+ix)\G(\g+ix)|^2$ ($\a,\g>0$ or $\bar\a=\g$ and $\Re(\a)>0$)
and are nowadays known as the symmetric
continuous Hahn polynomials. Atakishiyev and Suslov \cite{\AtakS, \S 3} have
shown that this is not the end of the story and have introduced
the continuous Hahn polynomials which have one extra parameter,
see also Askey \cite{\Aske}. These polynomials are orthogonal with respect
to the weight function $|\G(\a+ix)\G(\g+ix)|^2$ with $\Re(\a),\Re(\g)>0$ and
$\Im(\a)=-\Im(\g)$.

The goal of this paper is to show that Bateman's approach can be used
to prove the orthogonality relations for the continuous Hahn polynomials
by only using the Jacobi polynomials and the Fourier transform in a similar
way as Bateman \cite{\Batedrie} did in order to prove
\thetag{\vglorthoBatemanpols}. The orthogonality relations for the
continuous Hahn polynomials are not new, but this paper shows that
Bateman, Pasternack and Hardy could have found these orthogonality relations
if they had pursued their approach somewhat further. Moreover, this point of
view on the relationship between Jacobi polynomials and continuous Hahn
polynomials gives an intrinsic explanation for the occurrence of the
Jacobi polynomials in Atakishiyev and Suslov's proof \cite{\AtakS, \S 3}
of the orthogonality relation for the continuous Hahn polynomials.

We do not exactly know why Bateman and Hardy did not proceed to find the
continuous Hahn polynomials as early as the 1930's and 1940's. The
following
explanation has been communicated to me by Richard Askey, whom I thank for
letting me reproduce his view on this matter here. As to Hardy, we know that
he kept the special functions, when needed, as simple as possible and that
he only used special functions when he had to. So going beyond the Legendre
polynomial was no option to Hardy. Bateman, as most of his contemporaries,
thought of hypergeometric series as a function of the power series
variable $z$, cf. the definition in \S 2. But the argument of the
Bateman-Pasternack polynomial occurs in one of the parameters of a
hypergeometric series, cf. remark~2.2(i), and not in the power series
variable as is the case for e.g. the Jacobi polynomials, and thus they were
not in the line of thought at that time. This also explains why Rice
\cite{\Rice} looked for an appropriate generalisation of the Bateman and
Pasternack polynomial by introducing a variable at the power series spot.
However, Rice did not obtain orthogonal polynomials in this way.

There are more orthogonal systems
involving orthogonal polynomials that are mapped onto each other by the
Fourier transform, or by another integral transform such as the Mellin and
Hankel transform. The best known
example of this are the Hermite functions, i.e. Hermite polynomials
multiplied by $e^{-x^2/2}$, which are eigenfunctions of the Fourier
transform. For more examples we refer to Koornwinder \cite{\Kooreen},
\cite{\Koortwee} and to the integral transforms of the Bateman project
\cite{\ErdeTIT}. It is however important to note that in the
derivation presented here we do not use orthogonal systems, but
biorthogonal systems involving Jacobi polynomials. This gives the
possibility to introduce the necessary extra degrees of freedom.
Moreover, the result here is not a special case of Koornwinder's
result in which the Jacobi polynomials are mapped onto the Wilson
polynomials by use of the Jacobi function transform, cf.
\cite{\Kooreen}, \cite{\Koortwee}.

A striking aspect of \thetag{\vglBatemanpol} and
\thetag{\vglPasternackpol} is that the argument of the
orthogonal polynomial is a differential operator.
Badertscher and Koornwinder \cite{\BadeK}, see also \cite{\Koortwee},
have given group theoretic interpretations for several identities
involving orthogonal polynomials of differential operator argument.
In these cases these differential operators are acting on
spherical functions on Riemannian symmetric spaces, which are usually
more complicated special functions than just $\ch^{-1}x$.
A related paper in this direction is \cite{\Koel}.

The organisation of this paper is as follows. In \S 2 we derive the
Fourier transform of certain Jacobi polynomials in terms of continuous
Hahn polynomials and we discuss some applications. We also give the extension
of \thetag{\vglBatemanpol} and \thetag{\vglPasternackpol}, and we show how
some properties of the continuous Hahn polynomials can be derived from
properties of the Jacobi polynomials. In
\S 3 we prove the orthogonality relation for the
continuous Hahn polynomials from Parseval's identity for the Fourier
transform.

\demo{Acknowledgement} The work for this paper was initiated by the
referee report for \cite{\Koel}, in which the papers by Bateman and
Pasternack were mentioned. I thank this referee for drawing my attention
to these papers. \enddemo

\head\newsection. Fourier transform on Jacobi polynomials\endhead

The gamma function has been introduced by Euler in 1729 and is defined by
$$
\G (z) = \int_0^\infty t^{z-1} e^{-t}\, dt, \qquad \Re (z) >0.
$$
The fundamental recurrence relation $\G(z+1)=z\G(z)$ follows by
integration by parts.
A closely related integral is the beta integral;
$$
B (\a,\b) = \int_0^1 t^{\a-1} (1-t)^{\b-1} \, dt =
{{\G(\a)\G(\b)}\over{\G(\a+\b)}},
\qquad \Re(\a),\Re(\b)>0.
\tag\eqname{\vgldefbeta}
$$
The first proof of this result is given by Euler in 1772. More
information on proofs for the beta integral and related integrals and
sums as well as
references to the literature can be found in the papers
by Askey \cite{\Asketwee} and Rahman and Suslov \cite{\RahmS}.

The Jacobi polynomials $P_n^{(\a,\b)}(x)$ of degree $n$ in $x$
are the orthogonal polynomials with respect to
the beta measure shifted to the interval $[-1,1]$;
$$
\multline
\int_{-1}^1 (1-x)^\a (1+x)^\b P_n^{(\a,\b)}(x) P_m^{(\a,\b)}(x) \, dx =
\\
\d_{n,m} {{2^{\a+\b+1}}\over{2n+\a+\b+1}}
{{\G(n+\a+1)\G(n+\b+1)}\over{n!\, \G(n+\a+\b+1)}},
\endmultline
\tag\eqname{\vglorthoJacobi}
$$
for $\Re(\a)>-1$, $\Re(\b)>-1$. The orthogonality relations are ususally stated
for $\a>-1$, $\b>-1$, but they remain valid under this more general condition on
the parameters $\a$ and $\b$. The weight function is positive if and only if
$\a$ and $\b$ are real. An explicit expression for the Jacobi
polynomial $P_n^{(\a,\b)}(x)$ is given by a terminating hypergeometric series;
$$
P_n^{(\a,\b)}(x) = {{(\a+1)_n}\over{n!}} \, {}_2F_1\left(
{{-n,n+\a+\b+1}\atop{\a+1}}; {{1-x}\over 2}\right),
\tag\eqname{\vglJacobihyper}
$$
where the terminating hypergeometric series is defined by
$$
\gather
{}_{p+1}F_p\left( {{-n,a_1,\ldots,a_p}\atop{b_1,\ldots,b_p}};z\right) =
\sum_{k=0}^n {{(-n)_k(a_1)_k\ldots(a_p)_k}\over {(b_1)_k\ldots(b_p)_k}}
{{z^k}\over{k!}}, \qquad n\in\Zp, \\
(a)_k = a(a+1)(a+2)\ldots (a+k-1) = {{\G(a+k)}\over{\G(a)}}, \qquad k\in\Zp.
\tag\eqname{\vglPoch}
\endgather
$$
More information on Jacobi polynomials can be found in Szeg\H o's book
\cite{\Szeg, Ch.~4}.

Let us calculate a Fourier transform involving Jacobi polynomials.
Rewrite the Fourier transform
$$
\aligned
&\int_{-\infty}^\infty e^{-ixz} (1-\th x)^\a (1+\th x)^\b
P_n^{(\g,\d)}(\th x) \, dx \\
&= \int_{-1}^1 (1-t)^{\a-1+\hf iz} (1+t)^{\b-1-\hf iz} P_n^{(\g,\d)}(t)
\, dt \\
&= 2^{\a+\b-1} \int_0^1 u^{\a-1+\hf iz} (1-u)^{\b-1-\hf iz}
P_n^{(\g,\d)}(1-2u) \, du
\endaligned
\tag\eqname{\vglFJeen}
$$
using the substitutions $t=\th x$, $t=1-2u$. Use that
$$
{{dt}\over{dx}}=(\ch x)^{-2} = (1-\th x)(1+\th x),
\qquad
{{1-\th x}\over{1+\th x}} = e^{-2x}.
$$
Note that $x\mapsto (1-\th x)^\a (1+\th x)^\b P_n^{(\g,\d)}(\th x)\in
L^2(\R)\cap L^1(\R)$ for $\Re (\a), \Re(\b)> 0$.

In \thetag{\vglFJeen} we use the explicit series representation
\thetag{\vglJacobihyper} for
the Jacobi polynomial and the beta integral \thetag{\vgldefbeta}
to see that \thetag{\vglFJeen} equals
$$
\aligned
&2^{\a+\b-1} {{(\g+1)_n}\over{n!}} \sum_{k=0}^n
{{(-n)_k (n+\g+\d+1)_k}\over{k!\, (\g+1)_k}}
B(\a+k+\hf iz, \b-\hf iz) \\
&= 2^{\a+\b-1} {{(\g+1)_n}\over{n!}}
B(\a+\hf iz, \b-\hf iz)
\sum_{k=0}^n {{(-n)_k (n+\g+\d+1)_k (\a+\hf iz)_k}
\over{k!\, (\g+1)_k (\a+\b)_k}} \\
&= 2^{\a+\b-1} {{(\g+1)_n}\over{n!}}
B(\a+\hf iz, \b-\hf iz)
\, {}_3F_2 \left( {{-n,n+\g+\d+1,\a+\hf iz}\atop
{\g+1, \a+\b}};1\right).
\endaligned
\tag\eqname{\vglFJtwee}
$$
The identity obtained in this way can also be found in Erd\'elyi et al.
\cite{\ErdeTIT, 16.4(3)}.

The ${}_3F_2$-series in \thetag{\vglFJtwee}
is a continuous Hahn polynomial defined by, cf. \cite{\Aske},
$$
p_n(x;a,b,c,d) = i^n {{(a+c)_n(a+d)_n}\over{n!}}
\, {}_3F_2 \left( {{-n, n+a+b+c+d-1, a+ix}\atop{a+c,\ a+d}};1\right).
\tag\eqname{\vgldefcontHahn}
$$
So we have proved the following lemma.

\proclaim{Lemma \theoremname{\LemmaFJ}} For $z\in\R$, $\Re (\a), \Re(\b)> 0$
and $-\g\notin\N$ we have
$$
\gather
\int_{-\infty}^\infty e^{-ixz} (1-\th x)^\a (1+\th x)^\b
P_n^{(\g,\d)}(\th x) \, dx \\
= 2^{\a+\b-1} {{\G(\a+\hf iz)\G(\b-\hf iz)}\over{\G(\a+\b+n)}}
i^{-n} p_n(z/2; \a,\d-\b+1,\g-\a+1,\b),
\endgather
$$
where $P_n^{(\g,\d)}(x)$ is a Jacobi polynomial defined by
\thetag{\vglJacobihyper}
and $p_n(x;a,b,c,d)$ is a continuous Hahn polynomial defined by
\thetag{\vgldefcontHahn}.
\endproclaim

\demo{Remark \theoremname{\remFJ2}} (i) An equivalent
way of proving the lemma is by establishing
$$
\multline
p_n( -{i\over 2}\dx ;\a,\d-\b+1,\g-\a+1,\b)
\{ (1-\th x)^\a (1+\th x)^\b\} = \\
i^n (\a+\b)_n \{ (1-\th x)^\a (1+\th x)^\b\} P_n^{(\g,\d)}(\th x)
\endmultline
\tag\eqname{\vglcHdiffopJac}
$$
and applying the Fourier transform $\F$ to it and using
$\F( f^\prime)(z) = iz\F(f)(z)$.
Equation \thetag{\vglcHdiffopJac} can be proved from
$$
\gathered
\dx [ (1-\th x)^\a (1+\th x)^\b] = \bigl(\a+\b +(\a-\b)\th x\bigr)
[(1-\th x)^\a (1+\th x)^\b], \\
(\a+\hf \dx)_r [(1-\th x)^\a (1+\th x)^\b] =\qquad\qquad\qquad\qquad
\qquad\qquad\qquad\qquad\qquad
\\
\qquad\qquad\qquad\qquad\qquad\qquad\qquad
2^{-r} (1-\th x)^r (\a+\b)_r [(1-\th x)^\a (1+\th x)^\b],
\endgathered
\tag\eqname{\vgltweecHdiffopJac}
$$
which in turn can be proved by induction with respect to $r\in\Zp$.
Special cases of \thetag{\vglcHdiffopJac} are
\thetag{\vglBatemanpol} for $\a=\b={1\over 2}$, $\g=\d=0$ and
\thetag{\vglPasternackpol} for $\a=\b={1\over 2}(m+1)$, $\g=\d=0$ for
the Bateman polynomial and its generalisation by Pasternack. So we have
$$
\align
F_n^m(x) &= {1\over{i^n (m+1)_n}} p_n\bigl( -{i\over 2}x; \hf (1+m),
\hf (1-m),\hf (1-m),\hf (1+m)\bigr) \\
&= \, {}_3F_2\left( {{-n,n+1,\hf(1+m+x)}\atop{1,\ m+1}};1\right)
\endalign
$$
for the Bateman $(m=0)$ and Pasternack polynomials defined in
\thetag{\vglBatemanpol} and \thetag{\vglPasternackpol}.
The case $\a=\b={1\over 4}+{i\over 2}\lambda$, $\g=\d=-{1\over 2}$
of \thetag{\vglcHdiffopJac} has
been observed by Koornwinder \cite{\Koortwee, \S 2}.

(ii) Applying the inverse Fourier transform and taking $n=0$
gives the Fourier transform of $z\mapsto\G(\a+iz)\G(\b-iz)$, which is
closely related to the orthogonality measure for the Meixner-Pollaczek
polynomials, cf. Askey \cite{\Asketwee}. \enddemo

Lemma \LemmaFJ\ can be used to obtain identities for the continuous Hahn
polynomials from identities satisfied by the Jacobi polynomials. As a
first example we start with
the following generating functions for the Jacobi
polynomials, cf. \cite{\KoekS, (1.8.7), (1.8.6)},
$$
\gather
(1-t)^{-\g-\d-1} \, {}_2F_1\left(
{{{1\over 2}(\g+\d+1),{1\over 2}(\g+\d+2)}\atop{\g+1}};
{{2(x-1)t}\over{(1-t)^2}}\right) = \\ \sum_{n=0}^\infty
 {{(\g+\d+1)_n}\over{(\g+1)_n}} P_n^{(\g,\d)}(x) t^n,
\endgather
$$
and
$$
{}_0F_1\left( {{-}\atop{\g+1}};{{(x-1)t}\over 2}\right)\
{}_0F_1\left( {{-}\atop{\d+1}};{{(x+1)t}\over 2}\right) =
\sum_{n=0}^\infty {{P_n^{(\g,\d)}(x) t^n}\over{(\g+1)_n(\d+1)_n}}.
$$
A straightforward manipulation using lemma \LemmaFJ\ proves the
following generating functions for the continuous Hahn polynomials. The
first of these generating functions is also contained in \cite{\KoekS,
(1.4.8)}. We get
$$
\gather
(1-t)^{-\a-\b-\g-\d-1} \, {}_3F_2\left(
{{{1\over 2}(\a+\b+\g+\d-1),{1\over 2}(\a+\b+\g+\d),\a+iz}\atop
{\g+\a,\quad \a+\b}};
{{-4t}\over{(1-t)^2}}\right) = \\ \sum_{n=0}^\infty
 {{(\a+\b+\g+\d-1)_n}\over{(\a+\b)_n(\a+\g)_n}} (t/i)^n
p_n(z;\a,\d,\g,\b) ,
\endgather
$$
which is the generating function used by Bateman
\cite{\Bateeen, \S 3} and Pasternack \cite{\Past, (2.2)}, and
$$
\sum_{n=0}^\infty {{(t/i)^n p_n(z;\a,\d,\g,\b)}\over{
(\g+\a)_n(\d+\b)_n(\a+\b)_n}} =
\sum_{p=0}^\infty \sum_{k=0}^\infty {{(-t)^p t^k (\a+iz)_p
(\b-iz)_k}\over{p!\, (\g+\a)_p k!\, (\d+\b)_k (\a+\b)_{p+k}}}.
$$
The last series can be rewritten as a hypergeometric series in two variables
with arguments $x=-t$, $y=t$, cf. Appell and Kamp\'e de
F\'eriet \cite{\AppeKF, Ch.~IX, p.~150, (29)}.

As a second example we derive two relations between three continuous
Hahn polynomials from identities for Jacobi polynomials. First use
$\F( f^\prime)(z) = iz\F(f)(z)$, \thetag{\vgltweecHdiffopJac} and
$\dx P^{(\g,\d)}(x) = {1\over 2}(n+\g+\d+1) P_{n-1}^{(\g+1,\d+1)}(x)$,
cf. \cite{\Szeg, (4.21.7)},
and straightforward calculations to get after some rewriting
$$
\multline
(\a+\b+n) iz p_n(z;\a,\d,\g,\b) = (\a+\b) (\a+iz) p_n(z;\a+1,\d,\g-1,\b)\\
+ i(n+\a+\b+\g+\d-1)(\a+iz)(\b-iz) p_{n-1}(z;\a+1,\d,\g,\b+1).
\endmultline
$$
Another classical identity for the Jacobi polynomials, cf.
\cite{\Szeg, (4.5.4)},
$$
(n+\g+1) P_n^{(\g,\d)}(x) - (n+1) P_{n+1}^{(\g,\d)}(x) =
{1\over 2}(2n+\g+\d+2) (1-x) P_n^{(\g+1,\d)}(x)
$$
leads to
$$
\multline
(2n+\a+\b+\g+\d)(\a+iz) p_n(z;\a+1,\d,\g,\b) =\\
(\a+\b+n)(n+\g+\a) p_n(z;\a,\d,\g,\b) + i(n+1) p_{n+1}(z;\a,\d,\g,\b).
\endmultline
$$

It is also possible to use two or more identities for the Jacobi polynomials
in order to obtain identities for continuous Hahn polynomials. As an example
we indicate how the three-term recurrence relation for the continuous Hahn
polynomials can be derived, cf. Pasternack \cite{\Past, \S 5}.
Let $p_n$, cf. \thetag{\vglcHdiffopJac}, be defined by
$$
p_n\bigl(\dx\bigr)
\{ (1-\th x)^\a (1+\th x)^\b\} =
\{ (1-\th x)^\a (1+\th x)^\b\} P_n^{(\g,\d)}(\th x).
\tag\eqname{\vglcHdiffJactwee}
$$
Differentiate this identity once more to get
$$
\multline
\dx p_n\bigl(\dx\bigr)
\{ (1-\th x)^\a (1+\th x)^\b\} =
\{ (1-\th x)^\a (1+\th x)^\b\} \\
\times \bigl[ (\a+\b + (\a-\b)\th x) P_n^{(\g,\d)}(\th x)
+ (1-\th^2x) {{dP_n^{(\g,\d)}}\over{dx}}(\th x) \bigr].
\endmultline
$$
In the term in square brackets we use
$$
(1-x^2){{dP_n^{(\g,\d)}}\over{dx}}(x) =
A_n P^{(\g,\d)}_{n+1}(x) +
B_n P^{(\g,\d)}_n(x) +
C_n P^{(\g,\d)}_{n-1}(x),
$$
cf. Szeg\H o \cite{\Szeg, (4.5.5)} for the explicit values of the constants,
and the three-term recurrence relation for the Jacobi polynomials, cf.
\cite{\Szeg, (4.5.1)}, to get only Jacobi polynomials of degree $n+1$, $n$
and $n-1$. Recalling \thetag{\vglcHdiffJactwee} we find the three-term
recurrence relation for the continuous Hahn polynomials $p_n$. For the
explicit values of the coefficients we refer to \cite{\KoekS, (1.4.3)}.

\head\newsection. Orthogonality for the continuous Hahn polynomials
\endhead

The set of functions $x\mapsto (1-\th x)^\a (1+\th x)^\b
P_n^{(2\a+1,2\b+1)}$ is an orthogonal basis of $L^2(\R)$ and by
lemma \LemmaFJ\ it is mapped by the Fourier transform onto the set of
functions $z\mapsto \G(\a+\hf iz)\G(\b-\hf iz) p_n(\hf z;\a,\b,\a,\b)$.
Since the Fourier transform is isometric we obtain the orthogonality
relations for the continuous symmetric Hahn polynomials, cf.
\cite{\AskeWeen}, see also \cite{\Kooreen}, \cite{\Koortwee},
but we can do better as follows.

From the Parseval identity $2\pi\int_{\R} f(x)\bar g(x)dx =
\int_{\R}\bigl(\F f\bigr)(z) \overline{\bigl(\F g\bigr)}(z) dz$ for the
the Fourier transform for $f,g\in L^2(\R)$ we obtain
$$
\aligned
&2\pi \int_{-\infty}^\infty (1-\th x)^{\a+ a} (1+\th x)^{\b+ b}
P_n^{(\g,\d)}(\th x) P_m^{( c, d)}(\th x)\, dx \\
&= i^{m-n}2^{\a + a+\b + b-2}
\int_{-\infty}^\infty {{\G(\a+\hf iz)\G(\b-\hf iz)}\over{\G(\a+\b+n)}}
{{\G( a-\hf iz)\G( b+\hf iz)}\over{\G( a+ b+m)}} \\
&\qquad\times\,
p_n(z/2; \a,\d-\b+1,\g-\a+1,\b) \overline{p_m(z/2; \bar a,\bar d-\bar
b+1,\bar c-\bar a+1,\bar b)} \, dz \endaligned
\tag\eqname{\vglParseval1}
$$
for $\Re(\a),\Re(\b),\Re(a),\Re(b)>0$. Next we restrict the parameters in
the left hand side of \thetag{\vglParseval1} such that we get the orthogonality
relations \thetag{\vglorthoJacobi} for the Jacobi polynomials.
So we take $\Re(\a + a)>0$, $\Re(\b+ b)>0$, which is already satisfied, and
$\g= c=\a+ a-1$, $\d= d=\b+ b-1$. For these choices we put again
$t=\th x$ to see that the left hand side
\thetag{\vglParseval1} equals zero for $n\not= m$. The square norm follows
from \thetag{\vglorthoJacobi} and so we get
$$
\aligned
\int_{-\infty}^\infty \G(\a+iz)\G(\b-iz)
\G(a-iz)\G(b+iz) p_n(z; \a,b,a,\b)
\overline{p_m(z; \bar a,\bar\b,\bar\a,\bar b)} \, dz \\
= 2\pi \d_{n,m}
{{\G(\a+\b+n)\G(a+b+n)\G(n+\a+a)\G(n+\b+b)}\over{n!\, (2n+\a+\b+a+b-1)
\G(n+\a+\b+a+b-1)}}
\endaligned
\tag\eqname{\vglParseval2}
$$
subject to the conditions $\Re(\a),\Re(\b),\Re(a),\Re(b)>0$, since
$\g=c=\a+a-1$, $\d=d=\b+b-1$. If we replace $p_m$ in \thetag{\vglParseval2}
by $lc(p_m)z^m$ then \thetag{\vglParseval2} remains valid for
$0\leq m\leq n$. Here $lc(p_m)$ denotes the leading coefficient of
$\overline{p_m(z; \bar a,\bar\b,\bar\a,\bar b)}$, and since this equals
the leading coefficient of $p_n(z; \a,b,a,\b)$, cf.
\thetag{\vgldefcontHahn}, we may replace
$\overline{p_m(z; \bar a,\bar\b,\bar\a,\bar b)}$ by $p_m(z; \a,b,a,\b)$.
So we have proved the orthogonality relations for the continuous Hahn
polynomials stated in the next proposition
from the Parseval identity for the Fourier transform and from the explicit
knowledge of the Fourier transform of the Jacobi polynomial described in
lemma \LemmaFJ.

\proclaim{Proposition \theoremname{\ProporthocontHahn}} The continuous
Hahn polynomials defined in \thetag{\vgldefcontHahn} satisfy
$$
\multline
{1\over{2\pi}}\int_{-\infty}^\infty \G(\a+iz)\G(\b-iz)
\G(a-iz)\G(b+iz) p_n(z; \a,b,a,\b)
p_m(z;\a,b,a,\b) \, dz\\
= \d_{n,m}
{{\G(\a+\b+n)\G(a+b+n)\G(n+\a+a)\G(n+\b+b)}\over{n!\, (2n+\a+\b+a+b-1)
\G(n+\a+\b+a+b-1)}}
\endmultline
\tag\eqname{\vglorthocontHahn}
$$
for $\Re(\a),\Re(\b),\Re(a),\Re(b)>0$.
\endproclaim

\demo{Remark \theoremname{\remorthocontHahn}} (i) The weight function is
positive for $\bar a=\a$, $\bar b=\b$, or for $\a=\bar\b$, $a=\bar b$, which
follows from the invariance of \thetag{\vglorthocontHahn} under
interchanging $\a$ and $b$ or $\b$ and $a$.

(ii) The case $n=m=0$ of \thetag{\vglorthocontHahn} is Barnes' first
lemma from 1908, see e.g. Bailey \cite{\Bail, 1.7} or
Whittaker and Watson \cite{\WhitW, \S 14.52} for
the original proof by Barnes.
An equivalent proof of Barnes' first lemma as
given here can be found in Titchmarsh's book \cite{\Titc, (7.8.3)},
where the
Mellin transform is used instead of the Fourier transform.
Proposition \ProporthocontHahn\ can be obtained in a similar way using
the Mellin transform if we use the Jacobi polynomials of argument
$(1-x)/(1+x)$. We obtain an orthogonal system on $[0,\infty)$.
So we start with the Mellin transform
$$
\multline
\int_0^\infty {{x^\a}\over{(1+x)^{\a+\b}}} P_n^{(\g,\d)}\left(
{{1-x}\over{1+x}}\right) x^{-i\lambda-1}\, dx = \\
{{\G(\a-i\lambda)\G(\b-i\lambda)}\over{\G(\a+\b+n)}} i^{-n}
p_n(-\lambda;\a,\d-\b+1,\g-\a+1,\b)
\endmultline
\tag\eqname{\vglMellinJacobi}
$$
and the Parseval formula for the Mellin transform gives proposition
\ProporthocontHahn. See Hardy \cite{\Hard, \S 8} for this approach to
Bateman's polynomial $F_n$.
In \cite{\Koordrie, prop.~3.1} Koornwinder shows that the Laguerre
polynomials are mapped onto the Meixner-Pollaczek polynomials by the
Mellin transform. The Mellin transform of the underlying measures is
given in Titchmarsh \cite{\Titc, (7.8.1)}. This is a limiting case of
\thetag{\vglMellinJacobi}. Replace in \thetag{\vglMellinJacobi}
$x$ by $x/\d$ and $\b$ by $\d\xi$ with $\Re(\xi)>0$, $|\Im(\xi)|<\pi$
and let $\d\to\infty$. Similarly we can obtain the analogue of lemma
\LemmaFJ\ with the Laguerre and Meixner-Pollaczek polynomials by
a suitable limit transition.

(iii) The result \thetag{\vglorthocontHahn} in this form has been
proved first by Askey \cite{\Aske} using Barnes' first lemma and the
Chu-Vandermonde summation formula for a terminating ${}_2F_1$-series of unit
argument. Before that Atakishiyev and Suslov \cite{\AtakS}, see also
\cite{\NikiSU, \S 3.10.3.2}, proved \thetag{\vglorthocontHahn} in case of
a positive weight function. The method employed by Atakishiyev and Suslov
uses Barnes' first lemma, which is rewritten in terms of the beta integral
so that the orthogonality relations of the Jacobi polynomials can be used.

(iv) Another proof of proposition \ProporthocontHahn\ using
symmetry in the parameters $a$, $b$, $\a$, $\b$
is given by Kalnins and Miller \cite{\KalnM, \S 3}. They also
give a proof of Barnes' first lemma in this way.

(v) The continuous Hahn polynomials are not on the top shelf of the
Askey-scheme of hypergeometric polynomials. The most general hypergeometric
orthogonal polynomials with a continuous weight function are the Wilson
polynomials, cf. Wilson \cite{\Wils}, which have four degrees of freedom.
The continuous Hahn polynomials can be obtained from the Wilson polynomials
by a suitable limit process, cf. e.g. \cite{\KoekS, \S 2.2}.
\enddemo

The orthogonality relations stated in \S 1 are special cases of
\thetag{\vglorthocontHahn}. In particular, the orthogonality
\thetag{\vglorthoPasternackpols} for the Pasternack polynomials follows
by taking $\a=\b=\hf(1+m)$, $a=b=\hf(1-m)$. This shows also that we have
a positive weight function for $-1<m<1$ or $m\in i\R$, which are equivalent
for $m$ and $-m$. To see this we have to use the reflection formula
$\G(z)\G(1-z)= \pi\sin^{-1}(\pi z)$ and some straightforward manipulations
on goniometric and hyperbolic functions. It should also be noted that taking
the same values for the parameters in \thetag{\vglParseval2} gives
Bateman's (bi)orthogonality relations \thetag{\vglbiorthoPasternackpols}.

%%%%%%%%%%%%%%%%%References%%%%%%%%%%%%%%%%%%%%%%%%%%%%%%%%%%%%%%%%%%%%%%

\enddocument